\documentclass[11pt]{amsart}
\usepackage{amsmath}
\usepackage{amsfonts}
\usepackage{amsthm}
\usepackage[active]{srcltx}

\newtheorem{theorem}{Theorem}[section]

\newtheorem{remark}{Remark}[section]
\newtheorem{definition}{Definition}[section]
\newtheorem{corollary}{Corollary}[section]
\newtheorem*{TDW}{Theorem DW}

\numberwithin{equation}{section}

\newtheorem*{Go2}{Theorem G2}

\newtheorem*{Wa2}{Theorem W2}
\newtheorem*{Sa}{Theorem S}
\newtheorem*{Zh}{Theorem Zh}

\newtheorem*{Bakh1}{Theorem B1}
\newtheorem*{Bakh2}{Theorem B2}

\begin{document}
\author{Ushangi Goginava and Artur Sahakian}
\title[convergence and Summability of Fourier series ]{Convergence and
Summability of Multiple Fourier series and generalized variation}
\address{U. Goginava, Institute of Mathematics, Faculty of Exact and Natural
Sciences, Tbilisi State University, Chavchavadze str. 1, Tbilisi 0128,
Georgia}
\email{zazagoginava@gmail.com}
\address{A. Sahakian, Yerevan State University, Faculty of Mathematics and
Mechanics, Alex Manoukian str. 1, Yerevan 0025, Armenia}
\email{sart@ysu.am}
\maketitle

\begin{abstract}
In this paper we present results on convergence and Ces\`{a}ro summability of Multiple Fourier series of functions of bounded generalized variation.
\end{abstract}

\footnotetext{%
\par
The reseach of U. Goginava was supported by Shota Rustaveli National Science
Foundation grant no.13/06 (Geometry of function spaces, interpolation and
embedding theorems).}

\section{Classes of Functions of two variables of Bounded Generalized
Variation}

In 1881 Jordan \cite{Jo} introduced a class of functions of bounded
variation and applied it to the theory of Fourier series. This
notion was generalized hereinafter by many authors (quadratic variation, $\Phi$%
-variation, $\Lambda$-variation ets., see \cite{Wi}-\cite{Ch}). In two
dimensional case the class BV of functions of bounded variation was
introduced by Hardy \cite{Ha}.

In this section we introduce several classes of bivariate functions of Bounded Generalized Variation  and compare them with the class $HBV$ (se Definition \ref{def1.1} below), which is important for the applications in Fourier analysis (see Theorem S in Section 2.).

Let $f(x,y),\ (x,y)\in\mathbb{R}^2$ be a real function of two variables of period $2\pi $ with respect to
each variable. Given intervals $I=(a,b)$, $J=(c,d)$ and points $x,y$ from $%
T:=[0,2\pi ]$ we denote
\begin{equation*}
f(I,y):=f(b,y)-f(a,y),\qquad f(x,J)=f(x,d)-f(x,c)
\end{equation*}
and
\begin{equation*}
f(I,J):=f(a,c)-f(a,d)-f(b,c)+f(b,d).
\end{equation*}
Let $E=\{I_{i}\}$ be a collection of nonoverlapping intervals from $T$
ordered in arbitrary way and let $\Omega $ be the set of all such
collections $E$. Denote by $\Omega _{n}$ set of all collections of $n$
nonoverlapping intervals $I_{k}\subset T.$

For the sequence of positive numbers $\Lambda =\{\lambda
_{n}\}_{n=1}^{\infty }$ we define
\begin{equation*}
\Lambda V_{1}(f)=\sup_{y}\sup_{E\in \Omega }\sum_{n}\frac{|f(I_{i},y)|}{%
\lambda _{i}}\,\,\,\,\,\,\left( E=\{I_{i}\}\right) ,
\end{equation*}
\begin{equation*}
\Lambda V_{2}(f)=\sup_{x}\sup_{F\in \Omega }\sum_{m}\frac{|f(x,J_{j})|}{%
\lambda _{j}}\qquad (F=\{J_{j}\}),
\end{equation*}
\begin{equation*}
\Lambda V_{1,2}(f)=\sup_{F,\,E\in \Omega }\sum_{i}\sum_{j}\frac{%
|f(I_{i},J_{j})|}{\lambda _{i}\lambda _{j}}.
\end{equation*}

\begin{definition}\label {def1.1}
We say that the function $f$ has Bounded $\Lambda $-variation on $%
T^{2}=[0,2\pi ]^{2}$ and write $f\in \Lambda BV$, if
\begin{equation*}
\Lambda V(f):=\Lambda V_{1}(f)+\Lambda V_{2}(f)+\Lambda V_{1,2}(f)<\infty .
\end{equation*}%
We say that  $f$ has Bounded Partial $\Lambda $-variation and
write $f\in P\Lambda BV$ if
\begin{equation*}
P\Lambda V(f):=\Lambda V_{1}(f)+\Lambda V_{2}(f)<\infty .
\end{equation*}
\end{definition}

If $\lambda_n\equiv 1$ (or if $0<c<\lambda_n<C<\infty,\ n=1,2,\ldots$) the
classes $\Lambda BV$ and $P\Lambda BV$ coincide with the Hardy class $BV$
and PBV respectively. Hence it is reasonable to assume that $%
\lambda_n\to\infty$ and since the intervals in $E=\{I_i\}$ are ordered
arbitrarily, we will suppose, without loss of generality, that the sequence $%
\{\lambda_n\}$ is increasing. Thus,
\begin{equation}  \label{Lambda}
1<\lambda_1\leq \lambda_2\leq\ldots,\qquad \lim_{n\to\infty}\lambda_n=\infty.
\end{equation}

In the case when $\lambda _{n}=n,\ n=1,2\ldots $ we say \textit{Harmonic
Variation} instead of $\Lambda $-variation and write $H$ instead of $\Lambda$
($HBV$, $PHBV$, $HV(f)$, ets).

The notion of $\Lambda$-variation was introduced by D. Waterman \cite{W} in
one dimensional case and A. Sahakian \cite{Saha} in two dimensional case. The class $PBV$ as well as the class $PBV_{p}$ (see Definition \ref{def1.2}) was introdused by U. Goginava in \cite{GoEJA}.

\begin{definition}\label{def1.2}
Let $\Phi $-be a strictly increasing continuous function on $[0,+\infty )$
with $\Phi \left( 0\right) =0$. We say that the function $f$ has Bounded
Partial $\Phi $-variation on $T^{2}$ and write $f\in PBV_\Phi$, if
\begin{equation*}
V_{\Phi }^{\left( 1\right) }\left( f\right)
:=\sup\limits_{y}\sup\limits_{\{I_{i}\}\in \Omega
_{n}}\sum\limits_{i=1}^{n}\Phi \left( |f\left( I_{i},y\right) |\right)
<\infty ,\quad n=1,2,...,
\end{equation*}
\begin{equation*}
V_{\Phi }^{\left( 2\right) }\left( f\right)
:=\sup\limits_{x}\sup\limits_{\{J_{j}\}\in \Omega
_{m}}\sum\limits_{j=1}^{m}\Phi \left( |f\left( x,J_{j}\right) |\right)
<\infty ,\quad m=1,2,....
\end{equation*}
In the case when $\Phi \left( u\right) =u^{p},\,p\geq 1$, we say that $f$ has Bounded Partial $p$-variation and write $f\in PBV_{p}$.
\end{definition}
In the following theorem the necessary and sufficient conditions are obtained for the inclusion $P\Lambda BV\subset HBV$.
\begin{theorem}[U. Goginava, A. Sahakian \cite{EJA}]\label{T1}
Let $\Lambda =\{\lambda _{n}\}$ with $\lambda _{n}=n\gamma_n$ and $\gamma_{n}\geq
\gamma _{n+1}>0,\ n=1,2,....\,\,\,$. \newline
\indent
1) If
\begin{equation}  \label{T1-1}
\sum_{n=1}^{\infty }\frac{\gamma _{n}}{n}<\infty,
\end{equation}
\hskip10mm then $P\Lambda BV\subset HBV$.\newline
\indent 2) If  $\gamma_n=O(\gamma_{n^{[1+\delta]}})$ for some $\delta>0$ and
\begin{equation*}
\sum_{n=1}^{\infty }\frac{\gamma _{n}}{n}=\infty,
\end{equation*}
\hskip10mm then $P\Lambda BV\not\subset HBV$.
\end{theorem}

\begin{corollary}\label{cor1.1}
$PBV\subset HBV$ and $PHBV\not\subset HBV$.
\end{corollary}
\begin{corollary}\label{cor1.2}
Let $\Phi $ and $\Psi $ are conjugate functions in the sense of Yung ($%
ab\leq \Phi (a)+\Psi (b)$) and let for some $\{\lambda _{n}\}$ satisfying  (\ref{Lambda}),
\begin{equation}
\sum_{n=1}^{\infty }\Psi \left( \frac{1}{\lambda _{n}}\right) <\infty .
\label{Phi}
\end{equation}%
Then $PBV_{\Phi }\subset HBV$. In particular, $PBV_{p}\subset HBV$ for any $%
p>1$.
\end{corollary}

\begin{definition}[U. Goginava  \protect\cite{GoEJA}]
The Partial Modulus of Variation of a function $f$ are the functions $%
v_{1}\left( n,f\right)$ and $v_{2}\left( m,f\right) $ defined by
\begin{eqnarray*}
&&
v_{1}\left( n,f\right) :=\sup\limits_{y}\sup\limits_{\{I_{i}\}\in \Omega
_{n}}\sum\limits_{i=1}^{n}\left| f\left( I_{i},y\right) \right| ,\quad
n=1,2,\ldots,
\\&&
v_{2}\left( m,f\right) :=\sup\limits_{x}\sup\limits_{\{J_{k}\}\in \Omega
_{m}}\sum\limits_{i=1}^{m}\left| f\left( x,J_{k}\right) \right| ,\quad
m=1,2,\ldots.
\end{eqnarray*}
\end{definition}
For functions of one variable the concept of the modulus of variation was
introduced by Chanturia \cite{Ch}.
\begin{theorem}[U. Goginava, A. Sahakian \cite{EJA}]\label{T2}
Let $f$ be such that
\begin{equation*}
\sum\limits_{n=1}^{\infty }\frac{\sqrt{v_{j}\left( n,f\right) }}{n^{3/2}}%
<\infty ,\quad j=1,2.
\end{equation*}%
Then $f\in HBV.$
\end{theorem}
Another class of functions of generalized bounded variation was introduced by
M. Dyachenko and D. Waterman in \cite{DW}. Denoting by $\Gamma $ the the set of finite
collections of nonoverlapping rectangles $A_{k}:=\left[ \alpha _{k},\beta
_{k}\right] \times \left[ \gamma _{k},\delta _{k}\right] \subset T^{2}$ they
define
\begin{equation*}
\Lambda ^{\ast }V\left( f\right) :=\sup_{\{A_{k}\}\in \Gamma }\sum\limits_{k}%
\frac{\left\vert f\left( A_{k}\right) \right\vert }{\lambda _{k}}.
\end{equation*}

\begin{definition}[M. Dyachenko, D. Waterman \cite{DW}]
 We say that $f\in \Lambda ^{\ast }BV $ if
\begin{equation*}
\Lambda V(f):=\Lambda V_{1}(f)+\Lambda V_{2}(f)+\Lambda ^{\ast }V\left(
f\right)<\infty .
\end{equation*}
\end{definition}

In \cite{GMJ} we introduced a new classes of functions of generalized
bounded variation and investigate the convergence of Fourier series of
function of that classes.
For the sequence $\Lambda =\{\lambda _{n}\}_{n=1}^{\infty }$ we define
\begin{equation*}
\Lambda ^{\#}V_{1}(f)=\sup_{\{y_{i}\}\subset T}\sup_{\{I_{i}\}\in \Omega
}\sum_{i}\frac{|f(I_{i},{y_{i}})|}{\lambda _{i}},
\end{equation*}%
\begin{equation*}
\Lambda ^{\#}V_{2}(f)=\sup_{\{x_{j}\}\subset T}\sup_{\{J_{j}\}\in \Omega
}\sum_{j}\frac{|f(x_{j},J_{j}|}{\lambda _{j}}.
\end{equation*}
\begin{definition}[U. Goginava, A. Sahakian \cite{EJA}]
We say that $f\in\Lambda ^{\#}BV$, if
\begin{equation*}
\Lambda ^{\#}V(f):=\Lambda ^{\#}V_{1}(f)+\Lambda ^{\#}V_{2}(f)<\infty .
\end{equation*}
\end{definition}
It is not hard to see, that
\begin{equation}
\Lambda ^{\ast }BV\subset \Lambda ^{\#}BV\subset P\Lambda BV.  \label{cond0}
\end{equation}%
Obviously,  the function $f(x,y)=\mathrm{sign}(x-y)$ belongs
to  $P\Lambda BV\setminus \Lambda ^{\#}BV$ for any $\Lambda $. On the
other hand, we have proved the following result.
\begin{theorem}[U. Goginava, A. Sahakian \cite{GMJ}]
If $\Lambda =\{\lambda _{n}\}$ and%
\begin{equation*}
{\limsup\limits_{n\rightarrow \infty }}\left( \sum\limits_{k=1}^{n^{2}}%
\frac{1}{\lambda _{k}}\right)  \left( \sum\limits_{k=1}^{n}\frac{1}{%
\lambda _{k}}\right)^{-1} =+\infty ,
\end{equation*}
then $\Lambda ^{\#}BV\setminus \Lambda ^{\ast }BV\neq \emptyset $.
\end{theorem}
In the next theorem we characterize sequences  $\Lambda =\{\lambda _{n}\}$ for which the inclusion
$\Lambda ^{\#}BV\subset HBV$ holds.
\begin{theorem}[U. Goginava, A. Sahakian \cite{GMJ}]\label{TT}
Let $\Lambda =\{\lambda _{n}\}$.\newline
\indent a) If
\begin{equation*}
{\limsup\limits_{n\rightarrow \infty }}\frac{\lambda _{n}\log n}{n}<\infty ,
\end{equation*}%
\hskip10mm then%
\begin{equation*}
\Lambda ^{\#}BV\subset HBV.
\end{equation*}%
\indent b) If $\frac{\lambda _{n}}{n}\downarrow 0$ and%
\begin{equation*}
{\limsup\limits_{n\rightarrow \infty }}\frac{\lambda _{n}\log n }{n}=+\infty ,
\end{equation*}%
\hskip10mm  then%
\begin{equation*}
\Lambda ^{\#}BV\not\subset HBV.
\end{equation*}%
\end{theorem}

\begin{definition}[U. Goginava, A. Sahakian \cite{GMJ}]
Let $\Phi $-be a strictly increasing continuous function on $[0,+\infty )$
with $\Phi \left( 0\right) =0$. We say that the function $f\in B^{\#}V_{\Phi
}\left( T^{2}\right) $, if
\begin{equation*}
V_{\Phi ,1}^{\#}\left( f\right) :=\sup\limits_{\{y_{i}\}\subset
T}\sup\limits_{\{I_{i}\}\in \Omega }\sum\limits_{i}\Phi \left( |f\left(
I_{i},y_{i}\right) |\right) <\infty ,
\end{equation*}%
and%
\begin{equation*}
V_{\Phi ,2}^{\#}\left( f\right) :=\sup\limits_{\{x_{j}\}\subset
T}\sup\limits_{\{J_{j}\}\in \Omega }\sum\limits_{j}\Phi \left( |f\left(
x_{j},J_{j}\right) |\right) <\infty .
\end{equation*}
\end{definition}

Next, we define
\begin{equation*}
v_{1}^{\#}\left( n,f\right)
:=\sup\limits_{\{y_{i}\}_{i=1}^{n}}\sup\limits_{\{I_{i}\}\in \Omega
_{n}}\sum\limits_{i=1}^{n}\left\vert f\left( I_{i},y_{i}\right) \right\vert
,\quad n=1,2,\ldots ,
\end{equation*}%
\begin{equation*}
v_{2}^{\#}\left( m,f\right)
:=\sup\limits_{\{x_{j}\}_{j=1}^{m}}\sup\limits_{\{J_{k}\}\in \Omega
_{m}}\sum\limits_{j=1}^{m}\left\vert f\left( x_{j},J_{j}\right) \right\vert
,\quad m=1,2,\ldots .
\end{equation*}

\begin{theorem}[U. Goginava, A. Sahakian \cite{GMJ}]\label{phi}
Let $\Phi $ and $\Psi $ are conjugate functions in the sense of
Yung ($ab\leq \Phi (a)+\Psi (b)$) and let
\begin{equation*}
\sum_{n=1}^{\infty }\Psi\left(  \frac{\log n}n\right)
<\infty .
\end{equation*}%
Then%
\begin{equation*}
B^{\#}V_{\Phi }\subset \left\{ \frac{n}{\log n }\right\}
^{\#}BV.
\end{equation*}
\end{theorem}

\begin{theorem}[U. Goginava, A. Sahakian \cite{GMJ}]
\label{v(n)}Let%
\begin{equation*}
\sum\limits_{n=1}^{\infty }\frac{v_{s}^{\#}\left( f,n\right) \log n }{n^{2}}<\infty ,\qquad \ s=1,2.
\end{equation*}%
Then
$$f\in \left\{ \frac{n}{\log n }\right\} ^{\#}BV.$$
\end{theorem}
Observe that by Theorem \ref{TT} we have the inclusion $\left\{ \frac{n}{\log n }\right\} ^{\#}BV\subset HBV$.
Now, for a sequence $\Lambda=\{\lambda_n\}$ we denote
$$
\Lambda _{n}:=\left\{ \lambda _{k}\right\} _{k=n}^{\infty },\qquad n=1,2,\ldots
$$


\begin{definition}[U Goginava \cite{SMJ}]
We say that the function $f\in \Lambda^{\#} BV$\ is continuous in $\Lambda ^{\#}$-variation
and write $f\in C\Lambda ^{\#}V$, if%
\begin{equation*}
\lim\limits_{n\rightarrow \infty }\Lambda _{n}^{\#}V_{1}\left( f\right)
=\lim\limits_{n\rightarrow \infty }\Lambda _{n}^{\#}V_{2}\left( f\right) =0.
\end{equation*}%
\end{definition}

\begin{theorem}[U. Goginava, A. Sahakian \cite {CMA1}]
\label{continous} Let the sequence $\Lambda =\left\{ \lambda _{n}\right\} $
be such that
\begin{equation*}\label{maincond}
\liminf_{n\rightarrow \infty }\frac{\lambda _{2n}}{\lambda _{n}}=q>1.
\end{equation*}%
Then $\Lambda ^{\#}BV=C\Lambda ^{\#}V$.
\end{theorem}

\begin{theorem}[U. Goginava \cite{SMJ}]
\label{MV}Let $\alpha +\beta <1,\alpha ,\beta >0$ and%
\begin{equation*}
\sum\limits_{j=1}^{\infty }\frac{v_{s}^{\#}\left( f;2^{j}\right) }{%
2^{j\left( 1-\left( \alpha +\beta \right) \right) }}<\infty ,\qquad s=1,2.
\end{equation*}%
Then $f\in C\{n^{1-\left( \alpha +\beta \right) }\}^{\#}V$.
\end{theorem}

\section{Convergence of double Fourier series}
Everywhere in this and in the next section we suppose that the function $f$ is measurable on $\mathbb R^2$ and $2\pi $-periodic  with respect to each variable. The double Fourier series of a function  $f\in L^{1}\left( T^{2}\right)$ with respect to the trigonometric system is the series
\begin{equation*}
S \left[ f\right] :=\sum_{m,n=-\infty }^{+\infty }\widehat{f}\left(
m,n\right) e^{imx}e^{iny},
\end{equation*}
where
\begin{equation*}
\widehat{f}\left( m,n\right) =\frac{1}{4\pi ^{2}}\int_{0}^{2\pi
}\int_{0}^{2\pi }f(x,y)e^{-imx}e^{-iny}dxdy
\end{equation*}
are the Fourier coefficients of  $f$. The rectangular partial
sums of $S[f]$ are defined as follows:
\begin{equation*}
S_{M,N} \left[ f, (x,y)\right] :=\sum_{m=-M } ^{M } \sum_{n=-N }^N \widehat{f%
}\left( m,n\right) e^{imx}e^{iny},
\end{equation*}
In this paper we consider only {\bf  Pringsheim convergence } of double Fourier series, i.e {\bf convergence of rectangular partial
sums}   $S_{M,N} \left[ f, (x,y)\right]$, as $M,N\to\infty$.

We denote by $C(T^{2})$ the space of continuous on $\mathbb R^2$  and $2\pi $%
-periodic with respect to each variable  functions with the norm
\begin{equation*}
\Vert f\Vert _{C}:=\sup_{x,y\in T^{2}}|f(x,y)|.
\end{equation*}

For a function $f$ we denote by $f\left( x\pm 0,y\pm
0\right) $ the open coordinate quadrant limits (if exist) at the point $%
\left( x,y\right) $ and let $f^{\ast }(x,y)$
be the arithmetic mean of that quadrant limits:
\begin{multline}\label{limits}
f^{\ast }(x,y) :=\frac{1}{4}\left\{ f\left( x+0,y+0\right) +f\left(
x+0,y-0\right) \right.  \\
\left. +f\left( x-0,y+0\right) +f\left( x-0,y-0\right) \right\} .
\end{multline}
\begin{remark}
Observe that for a function $f\in \Lambda BV$ the quadrant limits $f\left( x\pm 0,y\pm 0\right) $
may not exist. As was shown in \cite{GMJ} for any function $f\in \Lambda ^{\#}BV$ the quadrant limits $f\left( x\pm 0,y\pm 0\right) $ exist at any point $\left( x,y\right) \in T^{2}$.
\end{remark}

We say the point $(x,y)\in T^2$ is a {\bf regular point} of a function $f$, if all quadrant limits in (\ref{limits}) exist.

The well known Dirichlet-Jordan theorem (see \cite{Zy}) states that the
Fourier series of a function $g(x), \ x\in T$ of bounded variation converges
at every point $x$ to the value $\left[ g\left( x+0\right)
+g\left(x-0\right) \right] /2$. If $g$ is in addition continuous on $T$ the
Fourier series converges uniformly on $T$.

Hardy \cite{Ha} generalized the Dirichlet-Jordan theorem to the double
Fourier series. He proved that if function $f$ has bounded variation in
the sense of Hardy ($f\in BV$), then $S \left[ f\right] $ converges  to
$f^{\ast }(x,y)$ at any regular point $\left( x,y\right) $.
If $f$ is in addition continuous on $T^{2}$ then $S
\left[ f\right] $ converges uniformly on $T^{2}$.

\begin{Sa} [Sahakian \cite{Saha}] \label{Sa}
{The Fourier series of a function $f\in HBV$ converges to $f^\ast(x,y)$
in any regular point $\left( x,y\right) $. The convergence is uniform on any
compact $K\subset T^{2}$, where the function $f$ is continuous.}
\end{Sa}
Theorem S was proved in \cite{Saha} under assumption that the
function is continuous on some open set containing $K$, while O. Sargsyan
noticed in \cite{SG}, that the continuity of $f$ on the compact $K$ is
sufficient.
\begin{definition}
We say that the class of functions $V \subset L^{1}(T^{2})$ is a class of
convergence on $T^{2}$, if for any function $f\in V $

 1) the Fourier series of $f$ converges to $f^{\ast }({x,y})$ at any
regular point $(x,y)$,

2) the convergence is uniform on any compact $K\subset T^{2}$,  where the function $f$ is continuous.
\end{definition}

The following results immediately follow from Theorems \ref{T1}, \ref{T2}, Corollary \ref{cor1.2} and Theorem S.

\begin{theorem}[U. Goginava, A. Sahakian \cite{EJA}]\label{d=2}
Let $\Lambda =\{\lambda _{n}\}$ with $\lambda _{n}=n\gamma_n$ and $\gamma_{n}\geq
\gamma _{n+1}>0,\ n=1,2,....\,\,\,$. \newline
\indent
1) If
\begin{equation*}
\sum_{n=1}^{\infty }\frac{\gamma _{n}}{n}<\infty,
\end{equation*}
then the class  $P\Lambda BV$ is a class of convergence on $T^{2}$.\newline
\indent 2) If $\gamma_n=O(\gamma_{n^{[1+\delta]}})$ for some $\delta>0$ and
\begin{equation*}  \label{T1-11}
\sum_{n=1}^{\infty }\frac{\gamma _{n}}{n}=\infty,
\end{equation*}
then
then there exists a continuous function $f\in P\Lambda BV$, the Fourier
series of which diverges over cubes at $\left( 0,0\right) .$
\end{theorem}
\begin{theorem}[U. Goginava, A. Sahakian \cite{EJA}]
The set of functions $f$ satisfying
\begin{equation*}
\sum\limits_{n=1}^{\infty }\frac{\sqrt{v_{j}\left( n,f\right) }}{n^{3/2}}%
<\infty ,\quad j=1,2,
\end{equation*}
is a class of convergence on $T^{2}$.
\end{theorem}

\begin{corollary}
The set of functions $f$ satisfying $v_{1}\left(n,f\right) =O\left(n^{\alpha }\right)$, $ v_{2}\left(n,f\right) =O\left(n^{\beta }\right),$ $0<\alpha ,\beta <1$,
is a class of convergence on $T^{2}$.
\end{corollary}

\begin{theorem} [U. Goginava \cite{GoEJA}]
The class  $PBV_{p},\ p\geq 1$, is a class of convergence on $T^{2}$.
\end{theorem}

From Theorem \ref{d=2} follows that for any $\delta >0$ the class
 $f\in P\left\{ \frac{n}{\log ^{1+\delta }n}\right\} BV$ is a class of
convergence. Moreover, one can not take here $\delta =0$. It is interesting to
compare this result with the following one obtained by M. Dyachenko and D. Waterman in
\cite{DW}.

\begin{TDW}[M. Dyachenko and D. Waterman \cite{DW}]
\label{TDW} If $f\in \left\{\frac n {\log n}\right\}^*BV$, then in any point
$(x,y)\in T^2$ the quadrant limits (\ref{limits}) exist and the double Fourier
series of $f$ converges to $f^\ast(x,y)$.\newline
Moreover, the sequence $\left\{\frac n {\log n}\right\}$ can not be replaced
with any sequence $\left\{\frac {n\alpha_n} {\log n}\right\}$, where $%
\alpha_n\to \infty$.
\end{TDW}

It is easy to show \textrm{(see\cite{DW})}, that $\left\{\frac n {\log
n}\right\}^*BV\subset HBV$, hence the convergence part of Theorem DW follows
from Theorem S. It is essential that the condition $f\in \left\{\frac n
{\log n}\right\}^*BV$ guaranties the existence of quadrant limits.

The following theorem immediately follows from Theorem \ref{TT}
and Theorem S.

\begin{theorem}[U. Goginava, A. Sahakian \cite{GMJ}]\label{GS}
If $\Lambda =\{\lambda _{n}\}$  and
\begin{equation*}
{\limsup\limits_{n\rightarrow \infty }}\frac{\lambda _{n}\log n}{n}<\infty ,  \label{cond1}
\end{equation*}%
then the class $\Lambda^{\#}BV$ is a class of convergence on $T^{2}$.

In particular, the class $\left\{ \frac{n}{\log n}\right\} ^{\#}BV$ is a class of convergence on $T^{2}$.
\end{theorem}

Theorem DW and (\ref{cond0}) imply that the sequence $\left\{ \frac{n}{\log
n}\right\} $ in Theorem \ref{GS} can not be replaced with any sequence $%
\left\{ \frac{n\alpha _{n}}{\log n}\right\} $, where $\alpha _{n}\rightarrow
\infty $.

Theorems \ref{phi}, \ref{v(n)} and \ref{GS} imply

\begin{theorem}[U. Goginava, A. Sahakian \cite{GMJ}]
The class $B^{\#}V_{\Phi }$ is a class of convergence on $T^{2}$, provided that (\ref{T1-1}) and (\ref{Phi}) hold.
\end{theorem}

\begin{theorem}[U. Goginava, A. Sahakian \cite{GMJ}]
Let%
\begin{equation*}
\sum\limits_{n=1}^{\infty }\frac{v_{s}^{\#}\left( f,n\right) \log n }{n^{2}}<\infty ,~\ \ s=1,2.
\end{equation*}%
Then in any point $(x,y)\in T^{2}$ the quadrant limits (\ref{limits}) exist
and the double Fourier series of $\ f$ converges to $f^\ast(x,y)$. The convergence is uniform on any compact $K\in
T^{2}$, if $\ f$ is continuous on $K$.
\end{theorem}

\section{Ces\`{a}ro
Summability of double Fourier
series}
For one-dimensional Fourier series D. Waterman has proved the following theorem.
\begin{Wa2}[D. Waterman \cite{Wa2}]\label{DW2}
Let $0<\alpha <1$. The Fourier series of a function $f\in \{n^{1-\alpha
}\}BV $ is everywhere $\left( C,-\alpha \right) $ bounded and is uniformly $%
\left( C,-\alpha \right) $ bounded on each closed interval of continuity of $%
f$. \newline
\indent If $f\in C\{n^{1-\alpha }\}BV$, then $S[f]$ is everywhere $\left(
C,-\alpha \right) $ summable to the value $\left[ f\left( x+0\right)
+f\left( x-0\right) \right] /2$ and the summability is uniform on each
closed interval of continuity.
\end{Wa2}

Later A.  Sablin proved in \cite{Sab}, that for $0<\alpha<1$ the classes $%
\{n^{1-\alpha }\}BV$ and $C\{n^{1-\alpha }\}BV$ coincide.

For double Fourier series  the {\bf Ces\`{a}ro  $(C;\alpha ,\beta )$-means} 
 of a function $f\in L^1(T^2)$ are defined by
\begin{equation*}
\sigma _{n,m}^{\alpha ,\beta }(f;x,y):=\frac{1}{A_{n}^{\alpha }}\frac{1}{%
A_{m}^{\beta }}\sum_{i=0}^{n}\sum_{j=0}^{m}A_{n-i}^{\alpha -1}A_{m-j}^{\beta
-1}S_{i,j}\left[ f,(x,y)\right],
\end{equation*}
where $\alpha ,\beta >-1$ and 
\begin{equation*}
A_{0}^{\alpha }=1,\,\,\,A_{k}^{\alpha }=\frac{(\alpha +1)\cdots (\alpha +k)}{%
k!},\quad k=1,2,....
\end{equation*}
The double Fourier series of $f$ is said to be $\left( C;\alpha
,\beta \right) $  {\bf summable to $s$ in a point }$\left( x,y\right) $, if
$$
\lim_{n,m\to \infty}\sigma _{n,m}^{\alpha ,\beta }(f;x,y)=s.
$$

L. Zhizhiashvili  has investigated the convergence of Ces\`{a}ro
means of double Fourier series of functions of bounded variation. In particular, the following
theorem was proved.

\begin{Zh}[L. Zhizhiashvili \protect\cite{Zh}]
If $f\in BV$, then the double Fourier series of $f$ is $\left( C;-\alpha
,-\beta \right) $ summable to $f^\ast(x,y)$
in any regular point $\left( x,y\right) $. The convergence is uniform on any
compact $K$, where the function $f$ is continuous.
\end{Zh}

For functions of partial bounded variation the problem was considered by the
first author.
\begin{Go2}[{U. Goginava} \cite{GoPRI}] Let $\alpha>0,\ \beta >0$.

1) If  $\alpha +\beta <1$, then for any $f\in C\left( T^{2}\right) \cap PBV$
the double Fourier series of  $f$
is uniformly $(C;-\alpha ,-\beta )$ summable to $f$.

2) If $\alpha +\beta \geq 1$, then there exists a
continuous function $f_{0}\in PBV$ such that the sequence $\sigma _{n,n}^{-\alpha ,-\beta }\left( f_{0};0,0\right) \,$ diverges.
\end{Go2}

In \cite{Szeg} we consider the following problem. \textit{Let }$\alpha
,\beta \in \left( 0,1\right) ,\,\alpha +\beta <1.$ \textit{Under what
conditions on the sequence }$\Lambda =\{\lambda _{n}\}$ \textit{the double
Fourier series of any function $f\in P\Lambda BV$ is $(C;-\alpha ,-\beta )$
summable?}

\begin{theorem}[U. Goginava, A. Sahakian \cite{Szeg}]
Let $\alpha ,\beta \in \left( 0,1\right) ,\ \alpha +\beta <1$ and let the sequence $\Lambda =\{\lambda _{k}\}$ be such that ${\lambda _{k}}{k^{\left( \alpha +\beta \right)-1 }}\downarrow 0$. 

1) If
\begin{equation*}
\qquad
\sum\limits_{k=1}^{\infty }\frac{\lambda _{k}}{k^{2-\left( \alpha +\beta
\right) }}<\infty,
\end{equation*}%
then the double Fourier series of any function $f\in P\Lambda BV$ is $\left(
C;-\alpha ,-\beta \right) $ summable to $f^*(x,y)$ at any regular point $\left( x,y\right) $. The
summability is uniform on any compact $K$, if  $f$ is continuous on the neighborhood of $K$.

2) If 
\begin{equation*}
\sum\limits_{k=1}^{\infty }\frac{\lambda _{k}}{k^{2-\left( \alpha +\beta
\right) }}=\infty,
\end{equation*}%
then there exists a continuous function $f\in P\Lambda BV$ for which the  $\left(
C;-\alpha ,-\beta \right) $ means of the double Fourier series diverges
over cubes at $\left( 0,0\right) .$
\end{theorem}

\begin{corollary}[U. Goginava, A. Sahakian \cite{Szeg}]
Let $\alpha ,\beta \in \left( 0,1\right) ,\,\alpha +\beta <1$. \newline

1) If $f\in P\left\{ \frac{n^{1-\left( \alpha +\beta \right) }}{\log
^{1+\varepsilon } n }\right\} BV$ for some $\varepsilon >0$,
then the double Fourier series of the function $f$ is $\left( C;-\alpha
,-\beta \right) $ summable to  $f^*(x,y)$
in any regular point $\left( x,y\right) $. The summability is uniform on any compact $K$, if
 $f$ is continuous on the neighborhood of $K$.

2) There exists a continuous function $f\in P\left\{ \frac{n^{1-\left(
\alpha +\beta \right) }}{\log \left( n+1\right) }\right\} BV$ such that $%
\left( C;-\alpha ,-\beta \right) $ means of two-dimensional Fourier series
of $f$ diverges over cubes at $\left( 0,0\right) .$
\end{corollary}

\begin{corollary}[U. Goginava, A. Sahakian \cite{Szeg}]
Let $\alpha,\beta \in \left( 0,1\right) ,\,\alpha +\beta <1$ \newline
and $\,\,f\in PBV$. Then the double Fourier series of the
function $f$ is $\left( C;-\alpha ,-\beta \right) $ summable to $f^*(x,y)$
 in any regular point $\left( x,y\right) $. The
summability is uniform on any compact $K$, if
 $f$ is continuous on the neighborhood of $K$.
\end{corollary}

In \cite{SMJ} the following problem was considred. \textit{Let }$\alpha ,\beta
\in \left( 0,1\right) ,\,\alpha +\beta <1.$ \textit{Under what conditions on
the sequence }$\Lambda =\{\lambda _{n}\}$ \textit{the double Fourier series
of any function $f\in C\Lambda ^{\#}BV$ is $(C;-\alpha ,-\beta )$ summable.}

\begin{theorem}[U. Goginava \cite{SMJ}]
\label{Main}a) Let $\alpha ,\beta \in \left( 0,1\right) ,\,\alpha +\beta <1$
and $f\in C\left\{ n^{1-\left( \alpha +\beta \right) }\right\} ^{\#}BV$.
Then the double Fourier series of  $f$ is $\left( C;-\alpha ,-\beta \right) $ summable to $f^*(x,y)$
 in any point $\left( x,y\right) $. The
summability is uniform on any compact $K\subset \mathbb{T}^{2}$, if
 $f$ is continuous on the neighborhood of $K$.\newline
b) Let $\Lambda :=\left\{ n^{1-\left( \alpha +\beta \right) }\xi
_{n}\right\} $, where $\xi _{n}\uparrow \infty $ as $n\rightarrow \infty $.
Then \ there exists a function $f\in C\left( \mathbb{T}^{2}\right) \cap
C\Lambda ^{\#}V$ for which $\left( C;-\alpha ,-\beta \right) $-means of
double Fourier series diverges unboundedly at $\left( 0,0\right) $.
\end{theorem}


Theorems \ref{continous}, \ref{MV} and  \ref{Main} imply the following results.

\begin{theorem}
Let $\alpha ,\beta \in \left( 0,1\right) ,\,\alpha +\beta <1$ and $f\in
\left\{ n^{1-\left( \alpha +\beta \right) }\right\} ^{\#}BV$. Then the double Fourier series of  $f$ is $\left( C;-\alpha ,-\beta \right) $ summable to $f^*(x,y)$
 in any point $\left( x,y\right) $.  The
summability is uniform on any compact $K\subset \mathbb{T}^{2}$, if
 $f$ is continuous on the neighborhood of $K$.\newline
\end{theorem}


\begin{theorem}
Let $\alpha ,\beta \in \left( 0,1\right) ,\,\alpha +\beta <1$ and%
\begin{equation*}
\sum\limits_{j=1}^{\infty }\frac{v_{s}^{\#}\left( f;2^{j}\right) }{%
2^{j\left( 1-\left( \alpha +\beta \right) \right) }}<\infty ,\quad s=1,2.
\end{equation*}%
Then the double Fourier series of  $f$ is $\left( C;-\alpha ,-\beta \right) $ summable to $f^*(x,y)$
 in any point $\left( x,y\right) $. The
summability is uniform on any compact $K\subset \mathbb{T}^{2}$, if $f$ is continuous on the neighborhood of $K$.\newline
\end{theorem}

\section{Classes of Functions of $d$ variables of Bounded Generalized
Variation}

Consider a function $f\left( x\right) $ defined on the
d-dimensional cube $T^{d}$ and a collection
of intervals
\begin{equation*}
J^{k}=\left( a^{k},b^{k}\right) \subset T,\qquad k=1,2,\ldots d.
\end{equation*}
For $d=1$ we set
\begin{equation*}
f\left( J^{1}\right) :=f\left( b^{1}\right) -f\left( a^{1}\right) .
\end{equation*}%
If for any function of $d-1$ variables the expression $f\left( J^{1}\times
\cdots \times J^{d-1}\right) $ is already defined, then for a function $f$ of $d$
variables the {\bf mixed difference} is defined as follows:
\begin{equation*}
f\left( J^{1}\times \cdots \times J^{d}\right) :=f\left( J^{1}\times \cdots
\times J^{d-1},b^{d}\right) -f\left( J^{1}\times \cdots \times
J^{d-1},a^{d}\right) .
\end{equation*}
For sequences of positive numbers
\begin{equation*}
\Lambda ^{j}=\{\lambda _{n}^{j}\}_{n=1}^{\infty },\quad
\lim_{n\to\infty}\lambda^j_n=\infty,\quad j=1,2,\ldots ,d,
\end{equation*}
and for a function $f(x)$, $x=(x_1,\ldots,x_d)\in T^d$ the $\left( \Lambda
^{1},\ldots ,\Lambda ^{d}\right) $-{\bf variation of $f$ with respect to
the index set }$D:=\{1,2,...,d\}$ is defined as follows:
\begin{equation*}
\left\{ \Lambda ^{1},\ldots ,\Lambda ^{d}\right\} V^{D}\left( f,T^{d}\right)
:=\sup\limits_{\{I_{i_{j}}^{j}\}\in \Omega }\ \sum\limits_{i_{1},...,i_{d}}%
\frac{\left\vert f\left( I_{i_{1}}^{1}\times \cdots \times
I_{i_{d}}^{d}\right) \right\vert }{\lambda^1 _{i_{1}}\cdots \lambda^d
_{i_{d}}}.
\end{equation*}

For an index set $\alpha =\{j_{1},...,j_{p}\}\subset D$ and any $x=\left(
x_{1},...,x_{d}\right) \in R^{d}$ we set ${\widetilde{\alpha }}:=D\setminus
\alpha $ and denote by $x_{\alpha }$ the vector of $R^{p}$ consisting of
components $x_{j},j\in \alpha $, i.e.
\begin{equation*}
x_{\alpha }=\left( x_{j_{1}},...,x_{j_{p}}\right) \in R^{p}.
\end{equation*}

By
\begin{equation*}
\left\{ \Lambda ^{j_{1}},...,\Lambda ^{j_{p}}\right\} V^{{\alpha }}\left(
f,x_{\widetilde{\alpha }},T^d\right) \quad \text {and}\quad f\left(
I_{i_{j_{1}}}^{1}\times \cdots \times I_{i_{j_{p}}}^{p},x_{\widetilde{\alpha
}}\right)
\end{equation*}
we denote respectively the $\left( \Lambda ^{j_{1}},...,\Lambda
^{j_{p}}\right) $-variation over the $p$-dimensional cube $T^{p}$ and mixed
difference of $f$ as a function of variables $x_{^{j_{1}}},...,x_{j_{p}}$
with fixed values $x_{^{\widetilde{\alpha }}}$ of other variables. The $%
\left( \Lambda ^{j_{1}},...,\Lambda ^{j_{p}}\right) ${\bf-variation of $%
f$\ with respect to the index set} ${\alpha }$ is defined as
follows:
\begin{equation*}
\left\{ \Lambda ^{j_{1}},...,\Lambda ^{j_{p}}\right\} V^{{\alpha }}\left(
f,T^{p}\right) =\sup\limits_{x_{^{{\widetilde{\alpha }}}}\in T^{d-p}}\left\{{%
\Lambda ^{j_{1}},...,\Lambda ^{j_{p}}}\right\}V^{{\alpha }}\left( f,x_{^{%
\widetilde{\alpha }}},T^{d}\right) .
\end{equation*}

\begin{definition}
We say that the function $f$ has total Bounded $\left( \Lambda
^{1},...,\Lambda ^{d}\right) $-variation on $T^{d}$ and write $f\in \left\{
\Lambda ^{1},...,\Lambda ^{d}\right\}BV\left( T^{d}\right) $, if
\begin{equation*}
\left\{ \Lambda ^{1},...,\Lambda ^{d}\right\}
V(f,T^{d}):=\sum\limits_{\alpha \subset D}\left\{ \Lambda ^{1},...,\Lambda
^{d}\right\}V^{{\alpha }}\left( f,T^{d}\right) <\infty .
\end{equation*}
\end{definition}

\begin{definition}
We say that the function $f$ is continuous in $\left( \Lambda
^{1},...,\Lambda ^{d}\right) $-variation on $T^{d}$ and write $f\in C\left\{
\Lambda ^{1},...,\Lambda ^{d}\right\} V\left( T^{d}\right) $, if%
\begin{equation*}
\lim\limits_{n\rightarrow \infty }\left\{ \Lambda ^{j_{1}},...,\Lambda
^{j_{k-1}},\Lambda _{n}^{j_{k}},\Lambda ^{j_{k+1}},...,\Lambda
^{j_{p}}\right\} V^{{\alpha }}\left( f,T^{d}\right) =0,\qquad k=1,2,\ldots ,p
\end{equation*}%
for any $\alpha \subset D,\ \alpha :=\{j_{1},...,j_{p}\}$, where $\Lambda
_{n}^{j_{k}}:=\left\{ \lambda _{s}^{j_{k}}\right\} _{s=n}^{\infty }$.
\end{definition}
The continuity of a function in $\Lambda$-variation was introduced by D. Waterman \cite{Wa2}
and was investigated in details by A. Bakhvalov (see \cite{Bakh}, \cite{Bakh2} and references therein).
This property is important for applications in the theory of Fourier series (see Theorem B1 in Section 5).
\begin{definition}
We say that the function $f$ has Bounded Partial $\left( \Lambda
^{1},...,\Lambda ^{d}\right) $-variation and write $f\in P\left\{ \Lambda
^{1},...,\Lambda ^{d}\right\} BV\left( T^{d}\right) $ if
\begin{equation*}
P\left\{ \Lambda ^{1},...,\Lambda ^{d}\right\}
V(f,T^{d}):=\sum\limits_{i=1}^{d}\Lambda ^{i}V^{\{i\}}\left( f,T^{d}\right)
<\infty .
\end{equation*}
\end{definition}
In the case when $\Lambda ^{1}=\cdots =\Lambda ^{d}=\Lambda $ we set
\begin{eqnarray*}
\Lambda BV(T^{d}):= &&\{\Lambda ^{1},...,\Lambda ^{d}\}BV(T^{d}), \\
C\Lambda V(T^{d}):= &&C\{\Lambda ^{1},...,\Lambda ^{d}\}V(T^{d}), \\
P\Lambda BV(T^{d}):= &&P\{\Lambda ^{1},...,\Lambda ^{d}\}BV(T^{d}).
\end{eqnarray*}%
If $\lambda _{n}=n$ for all $n=1,2\ldots $ we say \textit{Harmonic
Variation} instead of $\Lambda $-variation and write $H$ instead of $%
\Lambda, i.e. $ $HBV$, $PHBV$, $CHV$, ets.

\begin{theorem}[U. Goginava, A. Sahakian \cite{AM}] \label{41} 
Let $\Lambda =\{\lambda _{n}\}_{n=1}^{\infty }$ and $d\geq 2$. If
${\lambda _{n}}/{n}\downarrow 0$ and
\begin{equation*}\label{lambda}
\sum\limits_{n=1}^{\infty }\frac{\lambda _{n}\log ^{d-2}n}{n^{2}}<\infty ,
\end{equation*}%
then $P\Lambda BV(T^{d})\subset CHV(T^{d})$.
\end{theorem}
For a sequence $\Lambda =\{\lambda _{n}\}_{n=1}^{\infty }$ we denote%
\begin{equation*}
\Lambda ^{\#}V_{s}\left( f,T^{d}\right) :=\sup\limits_{\left\{ x^{i}{\left\{
s\right\} }\right\} \subset T^{d-1}}\sup\limits_{\left\{ I_{i}^{s}\right\}
\in \Omega }\sum\limits_{i}\frac{\left\vert f\left( I_{i}^{s},x^{i}{\left\{
s\right\} }\right) \right\vert }{\lambda _{i}},
\end{equation*}%
where
\begin{equation*}
x^{i}{\left\{ s\right\} }:=\left( x_{1}^{i},\ldots
,x_{s-1}^{i},x_{s+1}^{i},\ldots ,x_{d}^{i}\right) \quad \text{for}\quad
x^{i}:=\left( x_{1}^{i},\ldots ,x_{d}^{i}\right) .  \label{xis}
\end{equation*}

\begin{definition}
\label{def5} We say that  $f\in \Lambda
^{\#}BV\left( T^{d}\right) $, if%
\begin{equation*}
\Lambda ^{\#}V\left( f,T^{d}\right) :=\sum\limits_{s=1}^{d}\Lambda
^{\#}V_{s}\left( f,T^{d}\right) <\infty .
\end{equation*}
\end{definition}


\begin{theorem}[U. Goginava, A. Sahakian \cite{SBM}]
If $\Lambda =\left\{ \lambda _{n}\right\}$ with
\begin{equation*}\label{1}
\lambda _{n}=%
\frac{n}{\log ^{d-1}n },\quad n=2,3,\ldots,
\end{equation*}%
then $\Lambda
^{\#}BV(T^{d})\subset HBV(T^{d})$.
\end{theorem}

Now, we denote
\begin{equation*}
\Delta :=\{\delta =(\delta _{1},\ldots ,\delta _{d}):\delta _{i}=\pm 1,\
i=1,2,\ldots ,d\}  \label{Delta}
\end{equation*}%
and
\begin{equation*}
\pi _{\varepsilon \delta }(x):=(x_{1},\,x_{1}+\varepsilon \delta _{1})\times
\cdots \times (x_{d},\,x_{d}+\varepsilon \delta _{d}),
\end{equation*}%
for $x=(x_{1},\ldots ,x_{d})\in R^{d}$ and $\varepsilon >0$. We set $\pi
_{\delta }(x):=\pi _{\varepsilon \delta }(x)$, if $\varepsilon =1$.

For a function $f$ and $\delta\in \Delta$ we set
\begin{equation}  \label{lim}
f_\delta(x):=\lim_{t\in \pi_\delta(x),\ t\to x} f(t),
\end{equation}
if the last limit exists.

\begin{theorem}[U. Goginava, A. Sahakian \cite{SBM}]
\label{th2} Suppose  $\Lambda=\{\lambda_n\}$ and $f\in \Lambda ^{\#}BV\left( T^{d}\right)$.

a) If the limit $f_\delta(x)$ exists for some $x=(x_1,\ldots,x_d)\in T^d$
and some $\delta=(\delta_1,\ldots,\delta_d)\in \Delta$, then
\begin{equation*}  \label{th2a}
\lim\limits_{\varepsilon \rightarrow 0}\Lambda ^{\#}V\left(
f,\pi_{\varepsilon\delta}(x)\right) =0.
\end{equation*}

b) If $f$ is continuous on some compact $K\subset T^{d}$, then
\begin{equation*}
\lim_{\varepsilon \rightarrow 0}\Lambda ^{\#}V\left( f,\left[
x_{1}-\varepsilon ,x_{1}+\varepsilon \right] \times \cdots \times \left[
x_{d}-\varepsilon ,x_{d}+\varepsilon \right] \right) =0  \label{th2b}
\end{equation*}%
uniformly with respect to $x=(x_{1},\ldots ,x_{d})\in K$.
\end{theorem}

\begin{theorem}[U. Goginava, A. Sahakian \cite{SBM}]
If the function $f(x),\ x\in T^{d}$ satisfies the condition
\begin{equation*}
\sum\limits_{n=1}^{\infty }\frac{v_{s}^{\#}\left( f,n\right) \log
^{d-1}n }{n^{2}}<\infty ,\qquad \ s=1,2,...,d,
\end{equation*}%
then $f\in \left\{ \frac{n}{\log ^{d-1}n }\right\}
^{\#}BV\left( T^{d}\right) .$
\end{theorem}

\section{\protect\medskip Convergence of multiple Fourier series}

The Fourier series of function $f\in L^{1}\left( T^{d}\right) $ with respect
to the trigonometric system is the series
\begin{equation*}
S\left[ f\right] :=\sum_{n_{1},...,n_{d}=-\infty }^{+\infty }\widehat{f}%
\left( n_{1},....,n_{d}\right) e^{i\left( n_{1}x+\cdots +n_{d}x_{d}\right) },
\end{equation*}%
where
\begin{equation*}
\widehat{f}\left( n_{1},....,n_{d}\right) =\frac{1}{\left( 2\pi \right) ^{d}}%
\int_{T^{d}}f(x_{1},...,x_{d})e^{-i\left( n_{1}x_{1}+\cdots
+n_{d}x_{d}\right) }dx_{1}\cdots dx_{d}
\end{equation*}%
are the Fourier coefficients of $f$. The rectangular partial sums are
defined as follows:
$$
S_{N_{1},...,N_{d}}\left [f, (x_{1},...,x_{d})\right] =\sum_{n_{1}=-N_{1}}^{N_{1}}\cdots \sum_{n_{d}=-N_{d}}^{N_{d}}\widehat{f}%
\left( n_{1},....,n_{d}\right) e^{i\left( n_{1}x_{1}+\cdots
+n_{d}x_{d}\right) }
$$
We denote by $C(T^{d})$ the space of continuous and $2\pi $-periodic with
respect to each variable functions with the norm
\begin{equation*}
\Vert f\Vert _{C}:=\sup_{\left( x^{1},\ldots,\,x^{d}\right) \in
T^{d}}|f(x^{1},\ldots,x^{d})|.
\end{equation*}

We say that the point $x:=\left( x^{1},\ldots ,x^{d}\right) \in T^d$ is a {\bf
regular point} of a function $f$ if the limits (\ref{lim}) exist for all $\delta\in \Delta$.
For a regular point $x\in T^d$ we denote%
\begin{equation*}
f^{\ast }\left( x\right) :=\frac{1}{2^{d}}\sum_{\delta\in \Delta} f_\delta(x).  \label{limit}
\end{equation*}

\begin{definition}
We say that the class of functions $V \subset L^{1}(T^{d})$ is a class of
convergence on $T^{d}$, if for any function $f\in V $

1) the Fourier series of $f$ converges to $f^{\ast }({x})$ at any regular
point ${x}\in T^{d}$,

2) the convergence is uniform on any compact $K\subset T^{d}$, if $f$ is
continuous on the neighborhood of $K$.
\end{definition}
In \cite{Bakh} A. Bakhvalov showed that the class $HBV(T^{d})$ is not a class of
convergence on $T^{d}$, if $d>2$. On the other hand, he proved the following

\begin{Bakh1}[A. Bakhvalov \protect\cite{Bakh}]
\label{B} The class $CHV(T^{d})$ is a class of convergence on $T^{d}$ for any $d=1,2,\ldots $
\end{Bakh1}

Convergence of spherical and other partial sums of d-dimensional Fourier
series of functions of bounded $\Lambda $-variation was investigated in
deatails by M. Dyachenko \cite{D1,D2}, A. Bakhvalov \cite{Bakh,Bakh3}.

The first part of the next theorem is a consequence of Theorem \ref{41} and Theorem B1.

\begin{theorem}[[U. Goginava, A. Sahakian \cite{AM}]\label{t30}
Let $\Lambda =\{\lambda _{n}\}$ and $d\geq 2$.

a) If $\lambda_n/n \downarrow 0$ and
\begin{equation*}
\sum\limits_{n=1}^{\infty }\frac{\lambda _{n}\log ^{d-2}n}{n^{2}}<\infty ,
\end{equation*}%
then $P\Lambda BV$ is a class of convergence on $T^{d}$.

\medskip b) If
$ \frac{\lambda _{n}}{n}=O\left( \frac{\lambda _{\lbrack n^{\delta }]}}{%
[n^{\delta }]}\right)$
for some $\delta >1$, and
\begin{equation*}
\sum\limits_{n=1}^{\infty }\frac{\lambda _{n}\log ^{d-2}n}{n^{2}}=\infty ,
\end{equation*}%
then there exists a continuous function $f\in P\Lambda BV$, the Fourier
series of which diverges at $\left( 0,\ldots ,0\right) .$
\end{theorem}

Theorem \ref{t30} imply

\begin{corollary}
\label{c3} a) If $\Lambda =\left\{ \lambda _{n}\right\} _{n=1}^{\infty }$
with
\begin{equation*}
\lambda _{n}=\frac{n}{\log ^{d-1+\varepsilon }n},\qquad n=2,3,\ldots
\end{equation*}%
for some $\varepsilon >0$, then the class $P\Lambda BV$ is a class of
convergence on $T^{d}$.

\medskip b) If $\Lambda =\left\{ \lambda _{n}\right\} _{n=1}^{\infty }$ with
\begin{equation*}
\lambda _{n}=\frac{n}{\log ^{d-1}n},\qquad n=2,3,\ldots ,
\end{equation*}%
then the class $P\Lambda BV$ is not a class of convergence on $T^{d}$.
\end{corollary}

\begin{theorem}[Goginava, Sahakian \cite{SBM}]
\label{main} a) If $\Lambda =\left\{ \lambda _{n}\right\} _{n=1}^{\infty }$
with
\begin{equation*}
\lambda _{n}=\frac{n}{\log ^{d-1}n},\qquad n=2,3,\ldots,
\end{equation*}%
then the class $\Lambda ^{\#}BV\left( T^{d}\right) $ is a class of
convergence on $T^{d}$.

b)If $\Lambda =\left\{ \lambda _{n}\right\} _{n=1}^{\infty }$ with
\begin{equation*}
\lambda _{n}:=\left\{ \frac{n\xi _{n}}{\log ^{d-1} n}%
\right\} ,\quad n=2,3,\ldots,
\end{equation*}%
where $\xi _{n}\rightarrow \infty $ as $n\rightarrow \infty $, then there
exists a continuous function $f\in \Lambda ^{\#}BV\left( T^{d}\right) $ such
that the cubical partial sums of $d$-dimensional Fourier series of $f$
diverge unboundedly at $\left( 0,...,0\right) \in T^{d}$.
\end{theorem}

\begin{theorem}[Goginava, Sahakian \cite{SBM}]
For any $d>1$ the class of functions $f(x),\ x\in T^{d}$ satisfying the
following condition
\begin{equation*}
\sum\limits_{n=1}^{\infty }\frac{v_{s}^{\#}\left( f,n\right) \log
^{d-1}n }{n^{2}}<\infty ,~\ \ s=1,...,d,
\end{equation*}%
is a class of convergence.
\end{theorem}

\section{\protect Ces\`{a}ro summability of $d$-dimensional Fourier series}

The Ces\`{a}ro $\left( C;{\alpha }_{1},...,\alpha _{d}\right) \,$means of $d$%
-dimensional Fourier series of function $f\in L^1(T^d)$ is defined by
\begin{multline*}
\sigma _{m_{1},...,m_{d}}^{{\alpha }_{1},...,\alpha
_{d}}[f;(x_{1},...,x_{d})] \\
:=\left( \prod\limits_{i=1}^{d}A_{m_{i}}^{\alpha _{i}}\right)
^{-1}\sum\limits_{p_{1}=0}^{m_{1}}\cdots
\sum\limits_{p_{d}=0}^{m_{d}}\prod\limits_{i=1}^{d}A_{m_{i}-p_{i}}^{\alpha
_{i}-1}S_{p_{1},...,p_{d}}[f,(x_{1},...,x_{d})\mathbf{]},
\end{multline*}%
where{\large \ }
\begin{equation*}
A_{0}^{\alpha }=1,\,\,A_{n}^{\alpha }=\frac{\left( \alpha +1\right) \cdots
\left( \alpha +n\right) }{n!},\qquad \alpha>-1.
\end{equation*}
The Fourier series $S[f]$ is said to be
$\left( C;-\alpha _{1},...,-\alpha _{d}\right) $ summable to $s$
in a point $\left( x_1,\ldots,x_d\right)$, if
$$
\sigma _{m_{1},...,m_{d}}^{{\alpha }_{1},...,\alpha
_{d}}[f;(x_{1},...,x_{d})]\to s\quad\text{as}\quad x_{1},...,x_{d}\to \infty.
$$
\begin{definition}
We say that the class of functions $\Omega \subset L^{1}(T^{d})$ is a class
of $\left( C;-\alpha _{1},...,-\alpha _{d}\right) $ summability on $T^{d}$,
if the Cesaro $\left( C;{-\alpha }_{1},...,-\alpha _{d}\right) \,$ means of
Fourier series of any function $f\in \Omega $ converges to $f^{\ast }({x})$
at any regular point ${x}\in T^{d}$. The summability is uniform on any
compact $K\subset T^{d}$, if in addition, $f$ is continuous on the
neighborhood of $K$.
\end{definition}
The multivariate analog of Theorem W2 from Section 3 was proved by A. Bakhvalov in \cite{Bakh2}.
\begin{Bakh2}[A. Bakhvalov  \cite{Bakh2}] For any numbers $\alpha _{1},...,\alpha _{d}\in \left( 0,1\right)$ the class $C\{n^{1-\alpha_1}\},\ldots \{n^{1-\alpha_d}\}V(T^d)$ is a class of $\left( C;-\alpha _{1},...,-\alpha _{d}\right) $
summabi\-lity on $T^{d}$.
\end{Bakh2}

In the next theorem we consider the problem of  $\left( C;-\alpha _{1},...,-\alpha _{d}\right) $ summability of the Fourier series of functions of bounded partial $\Lambda$-variation.
\begin{theorem}[U. Goginava, A. Sahakian \cite{Stek}]\label{t3}
Suppose $\alpha _{1},...,\alpha _{d}\in \left( 0,1\right)$, $\alpha _{1}+\cdots
+\alpha _{d}<1$ and the sequence $\Lambda =\{\lambda _{n}\}_{n=1}^{\infty }$ is such that
\begin{equation*}
\frac{\lambda _{n}}{n^{1-\left( \alpha _{1}+\cdots +\alpha _{d}\right) }}%
\downarrow 0\text{ .}
\end{equation*}

a) If
\begin{equation*}
\sum\limits_{n=1}^{\infty }\frac{\lambda _{n}}{n^{2-\left( \alpha
_{1}+\cdots +\alpha _{d}\right) }}<\infty ,  \label{lambda1}
\end{equation*}%
then $P\Lambda BV(T^d)$ is a class of $\left( C;-\alpha _{1},...,-\alpha
_{d}\right) $ summability on $T^{d}$.

\medskip b) If
\begin{equation*}
\sum\limits_{n=1}^{\infty }\frac{\lambda _{n}}{n^{2-\left( \alpha
_{1}+\cdots +\alpha _{d}\right) }}=\infty ,  \label{lambda3}
\end{equation*}%
then there exists a continuous function $f\in P\Lambda BV(T^d)$ for which the sequence  $\sigma
_{N,...,N}^{-\alpha _{1},...,-\alpha _{d}}[f,\left( 0,...,0\right) ]$
diverges.
\end{theorem}

\begin{corollary}[U. Goginava, A. Sahakian \cite{Stek}]
Suppose $\alpha _{1},...,\alpha _{d}\in \left( 0,1\right),$ $\alpha _{1}+\cdots
+\alpha _{d}<1$ and $\Lambda =\left\{ \lambda _{n}\right\} _{n=1}^{\infty }$.

a) If
\begin{equation*}
\lambda _{n}=\frac{n^{1-\left( \alpha _{1}+\cdots +\alpha _{d}\right) }}{%
\log ^{1+\varepsilon }n },\qquad n=2.3.\ldots
\end{equation*}%
for some $\varepsilon >0$, then the class $P\Lambda BV(T^d)$ is a class of $%
\left( C;-\alpha _{1},...,-\alpha _{d}\right) $ summability on $T^{d}$.

b) If
\begin{equation*}
\lambda _{n}=\frac{n^{1-\left( \alpha _{1}+\cdots +\alpha _{d}\right) }}{%
\log n },\qquad n=2.3.\ldots,
\end{equation*}%
then $P\Lambda BV(T^d)$ is not a class of $\left( C;-\alpha
_{1},...,-\alpha _{d}\right) $ summability on $T^{d}$.
\end{corollary}

\begin{theorem}[U. Goginava, A. Sahakian  \cite{Stek}]
Let $\alpha _{1},...,\alpha _{d}\in \left( 0,1\right),$ $\alpha _{1}+\cdots
+\alpha _{d}<1.$ Then the set of functions $f$ satisfying the conditions
\begin{equation*}
\sum\limits_{j=0}^{\infty }\frac{\left( v_{i}\left(
2^{j},f\right) \right) ^{\alpha _{i}/\left( \alpha _{1}+\cdots +\alpha
_{d}\right) }}{2^{j\left( \alpha _{i}/\left( \alpha _{1}+\cdots +\alpha
_{d}\right) -\alpha _{i}\right) }}<\infty\quad \text{for}\quad i=1,...,d,
\end{equation*}%
is a class of $\left( C;-\alpha _{1},...,-\alpha _{d}\right) $ summability
on $T^{d}$.
\end{theorem}

\begin{theorem}[U. Goginava, A. Sahakian  \cite{Stek}]
Suppose  $\alpha _{1},...,\alpha _{d}\in \left( 0,1\right)$, $\alpha _{1}+\cdots
+\alpha _{d}<1/p$, $\, p\geq 1.$ Then the class $PBV_{p}$ is a class of $\left(
C;-\alpha _{1},...,-\alpha _{d}\right) $ summability on $T^{d}$.
\end{theorem}

In \cite{GoPRI} the first author has proved that the class $PBV_{p}$ is not a
class of $\left( C;-\alpha _{1},...,-\alpha _{d}\right) $ summability on $%
T^{d}$, if  $\alpha _{1},...,\alpha _{d}\in \left( 0,1\right),$ and $\alpha
_{1}+\cdots +\alpha _{d}\geq 1/p.$

\begin{corollary}[U. Goginava, A. Sahakian  \cite{Stek}]
Suppose  $\alpha _{1},...,\alpha _{d}\in \left( 0,1\right)$, $\alpha _{1}+\cdots
+\alpha _{d}<1$. Then the set of functions $f$ satisfying 
\begin{equation*}
v_{i}\left( 2^{j},f\right) =O\left( 2^{j\gamma }\right)\quad \text{for}\quad i=1,...,d,
,\qquad
\end{equation*}%
 is a class of $\left( C;-\alpha _{1},...,-\alpha _{d}\right) $ summability
on $T^{d}$.
\end{corollary}

\end{document}